\input amstex
\documentstyle{amsppt}
\NoRunningHeads
\pagewidth{14cm}
\pageheight{19.1cm}
\magnification=1200
\loadmsbm
\loadbold
\define\ch{\operatorname {ch.}}
\define\aut{\operatorname {Aut}}

\define\rank{\operatorname {rank}}

\define\bs{\boldsymbol}

\define\ol{\overline}
\define\vvp{\varphi}

\define\bk{\boldkey}
\define\lpt{[t_1^{\pm 1},\ldots, t_n^{\pm 1}]}

\NoBlackBoxes

\author 
Yoji Yoshii
\footnote{\leftline{Research supported by a PIMS Postdoctoral
Fellowship (2001).}}
\endauthor

\affil
Department of Mathematical Sciences, \\
University of Alberta, \\
Edmonton, Alberta, \\
T6G 2G1  Canada \\
yoshii\@math.ualberta.ca
\endaffil
\topmatter
\title
Classification of quantum tori with involution
\endtitle
\abstract
Quantum tori with graded involution
appear as coordinate algebras of extended affine Lie algebras
of type $\text A_1$, $\text C$ and $\text {BC}$.  We classify them in the category of
algebras with involution.
From this, we obtain
precise information on 
the root systems of 
extended affine Lie algebras of type $\text C$.
\endabstract
\endtopmatter

\document
\baselineskip 3.4ex

\head
Introduction
\endhead

Let $F$ be a field.
A quantum torus $F_{\bs q}$ is a noncommutative analogue 
of the algebra of Laurent polynomials over $F$, determined by a certain 
$n\times n$ matrix
$\bs q$.
Quantum tori appeared in several areas, e.g. 
quantum affine varieties [6], extended affine Lie algebras
[5]
or quantum physics
[7].
In noncommutative geometry or quantum physics, a special type of quantum tori called
a {\it noncommutative torus} is considered
(see Remark 0.1).

Our first purpose in this paper is to classify the graded involutions of
quantum tori.
It is known [1]
that the existence of a graded involution of $F_{\bs q}$
is equivalent to $\bs q$ being elementary,
i.e., all the entries of $\bs q$ are $1$ or $-1$.
We prove that for an elementary $\bs q$
we have
$F_{\bs q}\cong F_{\bs h_{l,n}}$, where
$$\bs h_{l,n}=
\overbrace
{\bs h\times\cdots\times\bs h}^{\text{$l$-times}}
\times 
\bk 1_{n-2l}
\quad\text{and}\quad
\bs h=
\left (\matrix \format\r&\quad\r \\
    1   & -1      \\
   -1   &  1          
\endmatrix
\right )
\quad\text{(Theorem 1.10)}
$$
(see Definition 1.4 for the notation $\times$).
Then we classify graded involutions $\tau$
of the elementary quantum torus $F_{\bs h_{l,n}}$.
We obtain that the algebra with involution
$(F_{\bs h_{l,n}},\tau)$ is isomorphic to
$$
(F_{\bs h_{l,n}},*),\quad 
(F_{\bs h_{l,n}},\tau_{1})\quad
\text{or}\quad
(F_{\bs h_{l,n}},\tau_{2})
\quad\text{(Theorem 2.7)}
$$
for three unique involutions
$*$, $\tau_1$ and $\tau_2$.

A quantum torus has a natural $\Bbb Z^n$-grading.
For any graded involution the subset of $\Bbb Z^n$,
consisting of the degrees in which homogeneous elements are fixed by the involution,
is a so-called semilattice, studied in [1].
In Lemma 4.1 we determine the index, an invariant of any semilattice [4],
for each of the 3 involutions of Theorem 2.7.
As a result,
the 3 semilattices are pairwise non-similar.
Moreover, 
we introduce a natural similarity invariant of semilattices called {\it saturation number}
(Definition 4.2).
Using this concept, we show 
that $l$ in the three semilattices above is a similarity invariant.
This allows us to complete the classification of semilattices
determined by quantum tori with graded involution
(Theorem 4.6).

Quantum tori with graded involution
appear as coordinate algebras of extended affine Lie algebras
of type $\text A_1$ in [11], $\text C$ in [2] 
and $\text {BC}$ in [3].  
Isomorphic coordinate algebras give rise to isomorphic
extended affine Lie algebras.
Thus, our results provide a finer classification of
extended affine Lie algebras in the above types.
Also, we obtain more precise information on 
the difference between
extended affine root systems
and 
the root systems of 
extended affine Lie algebras of type $\text C_r$ for $r\geq 3$
than the one described in [2]
(see Corollary 5.4 and 5.5).

The organization of the paper is as follows.
In \S 1 we define elementary quantum tori
and classify them.
In \S 2 we classify
(elementary) quantum tori with involution.
In \S 3 we review semilattices.
In \S 4 we obtain the classification of
semilattices determined by
(elementary) quantum tori with involution.
In the final section 
extended affine root systems of type $\text C$
are reviewed
and the difference to the root systems of extended affine Lie algebras
of type $\text C$ is discussed.

This is part of my Ph.D thesis, written at the University of Ottawa.
I would like to thank my supervisor, Professor Erhard Neher,
for his encouragement and suggestions.

\head
\S 1  Elementary quantum tori
\endhead

We begin by recalling quantum tori
(see [8]).
An $n\times n$
matrix ${\boldsymbol q}=(q_{ij})$ over a field $F$
such that  $q_{ii}=1$ and $q_{ji}=q_{ij}^{-1}$
is called a {\it quantum data matrix}
or simply a {\it quantum matrix}.
(This notion should not be confused with the use of the word
``quantum matrix'' in quantum algebra,
see e.g. [9].  But in our argument, no confusion will arise,
and so
we will simply call the $\bs q$ a quantum matrix.)
The {\it quantum torus 
$F_{\boldsymbol q}
=F_{\boldsymbol q}[t^{\pm 1}_{1},\ldots,t^{\pm 1}_{n}]$ 
determined by} a quantum matrix $\boldsymbol q$
is defined as
the associative algebra over $F$
with $2n$ generators $t^{\pm 1}_{1},\ldots,t^{\pm 1}_{n}$,
and relations
$t_{i}t_{i}^{-1}=t_{i}^{-1}t_{i}=1$ and
$t_{j}t_{i}=q_{ij}t_{i}t_{j}$ for all $1\leq i,j\leq n$.
Note that
$F_{\boldsymbol q}$ is commutative
if and only if $\boldsymbol q=\boldkey 1$ where
all the entries of $\boldkey{1}$ 
are $1$.
In this case, the quantum torus 
$F_{\bk 1}$
becomes
the algebra $F[t^{\pm 1}_{1},\ldots,t^{\pm 1}_{n}]$
of Laurent polynomials.
\remark{Remark 1.0}
For $F=\Bbb C$,
if we assume that $|q_{ij}|=1$ for all $i,j$,
then $\Bbb C_{\boldsymbol q}$ is a noncommutative torus 
[10].
Let $\theta_{ij}\in\Bbb R$ be such that $q_{ij}=e^{2\pi i\theta_{ij}}$.
Then $\bs\theta=(\theta_{ij})$ is an antisymmetric matrix over $\Bbb R$.
In noncommutative geometry or quantum physics,
one studies
the $C^*$-algebra completion of
the quantum torus as defined above
(see e.g. [10] or [7]).
\endremark

\enskip

Let $\Lambda=\Lambda_n$ be the free abelian group of rank $n$.
We give a $\Lambda$-grading of
the quantum torus
$F_{\bs q}=F_{\bs q}\lpt$
in the following way:
For any basis
$\{\bs\sigma_{1},\ldots,\bs\sigma_{n}\}$
of $\Lambda$,
we define the degree of
$$t_{\bs\alpha}:=t_1^{\alpha_1}\cdots t_n^{\alpha_n}
\quad\text{for $\bs\alpha
=\alpha_{1}\bs\sigma_{1}+\cdots+\alpha_{n}\bs\sigma_{n}\in\Lambda$
as $\bs\alpha$.}
$$
Then 
$F_{\bs q}=\oplus_{\bs\alpha\in\Lambda}\ Ft_{\bs\alpha}$ 
becomes a $\Lambda$-graded algebra.
We call this grading
the {\it toral $\Lambda$-grading of $F_{\bs q}$
determined by $\langle\bs\sigma_{1},\ldots,\bs\sigma_{n}\rangle$}.
Sometimes it is referred to as a 
{\it $\langle\bs\sigma_{1},\ldots,\bs\sigma_{n}\rangle$-grading}.
Also, if we write
$F_{\bs q}=\oplus_{\bs\alpha\in\Lambda}\ (F_{\bs q})_{\bs\alpha}$
or
$F_{\bs q}=\oplus_{\bs\alpha\in\Lambda}\ Ft_{\bs\alpha}$,
we are assuming some toral $\Lambda$-grading of $F_{\bs q}$.
One can check that the multiplication rule in
$F_{\boldsymbol q}$
for
$\bs q=(q_{ij})$
is the following:
for 
$\bs\beta
=\beta_{1}\bs\sigma_{1}+\cdots+\beta_{n}\bs\sigma_{n}\in\Lambda$,
$$t_{\bs\alpha}
t_{\bs\beta}
=\prod_{i<j}\ 
q_{ij}^{\alpha_{j}\beta_{i}}
t_{\bs\alpha+\bs\beta}. \tag1.1$$

\proclaim{Lemma 1.2}
If 
$\vvp \: F_{\bs q}=\oplus_{\bs\alpha\in\Lambda}\ Ft_{\bs\alpha}
\tilde{\longrightarrow} 
F_{\bs\eta}=\oplus_{\bs\alpha\in\Lambda}\ Ft_{\bs\alpha}$ 
is an isomorphism of algebras,
then there exists
the induced group automorphism $p$ of $\Lambda$
such that
$\vvp (Ft_{\bs\alpha})
=Ft_{p(\bs\alpha)}$ for all $\bs\alpha\in\Lambda$.
\endproclaim
\demo{Proof}
It is easily seen that
the units of any quantum torus with toral grading
are nonzero homogeneous elements.
Thus, since $\vvp (t_{\bs\alpha})$ is a unit
for any $\bs\alpha\in\Lambda$, there
exists $p(\bs\alpha)\in\Lambda$ such that
$\vvp (Ft_{\bs\alpha})=Ft_{p(\bs\alpha)}$,
and the map $p\:\Lambda\longrightarrow \Lambda$
is well-defined.
It is straightforward to check that
$p$ is an automorphism of $\Lambda$.
\qed
\enddemo

For quantum matrices
$\bs q$ and $\bs\eta$,
we say that 
{\it $\bs q$ is equivalent to $\bs\eta$}
and denote this by $\bs q\cong\bs\eta$
if $F_{\bs q}\cong F_{\bs\eta}$.
This is an equivalence relation.
Note that $\bs q\cong\bk 1$ implies $\bs q=\bk 1$.

If $F_{\bs q}$ has a toral $\Lambda$-grading, 
the centre $Z(F_{\bs q})$ of $F_{\bs q}$ is graded by some
subgroup of $\Lambda$ which we call the
{\it grading subgroup} of $Z(F_{\bs q})$.
 If $F_{\bs q}$ and $F_{\bs\eta}$ each have toral $\Lambda$-gradings,
we write $F_{\bs q}\cong_{\Lambda} F_{\bs\eta}$
to mean that $F_{\bs q}$ and $F_{\bs\eta}$ are isomorphic as
$\Lambda$-graded algebras.  Moreover, in that case
the grading subgroups of $Z(F_{\bs q})$ and $Z(F_{\bs\eta})$
coincide.

\proclaim{Lemma 1.3}
Let $\bs q$ and $\bs\eta=(\eta_{ij})_{1\leq i,j\leq n}$
be quantum matrices,
and let
$F_{\bs q}$ respectively $F_{\bs\eta}$
be the corresponding quantum tori.
Then
the following are equivalent:

\rom{(i)} $\bs q\cong\bs\eta$, i.e., 
$F_{\bs q}\cong F_{\bs\eta}$ as algebras,

\rom{(ii)} for any toral grading of $F_{\bs q}$,
there exists a basis
$\langle\bs\sigma_{1},\ldots,\bs\sigma_{n}\rangle$ 
of $\Lambda$ and
nonzero homogeneous elements $x_{i}
\in F_{\bs q}$ of degree $\bs\sigma_{i}$
such that
$x_{j}x_{i}=\eta_{ij}x_{i}x_{j}$ 
for all $1\leq i<j\leq n$,

\rom{(iii)} for any toral grading of $F_{\bs q}$,
there exists a toral grading
of $F_{\bs\eta}$
such that
$F_{\bs q}\cong_{\Lambda} F_{\bs\eta}$. 
In that case,
the grading subgroups of the centres $Z(F_{\bs q})$ and $Z(F_{\bs\eta})$
coincide.
\endproclaim
\demo{Proof}
We prove (i) $\Longrightarrow$ (ii) $\Longrightarrow$
(iii) $\Longrightarrow$ (i).
Suppose that (i) holds, i.e.,
there exists an isomorphism $\vvp$ from 
$F_{\bs q}$
onto $F_{\bs\eta}$.
Give a toral $\Lambda$-grading to
$F_{\bs q}$ 
so that 
$F_{\bs q}=\oplus_{\bs\alpha\in\Lambda}\ (F_{\bs q})_{\bs\alpha}$
and a toral $\langle\bs\varepsilon_{1},\ldots,\bs\varepsilon_{n}\rangle$-grading to
$F_{\bs\eta}=F_{\bs\eta}[t_1^{\pm 1},\ldots,t_n^{\pm 1}]$
so that
$F_{\bs\eta}=\oplus_{\bs\alpha\in\Lambda}\ Ft_{\bs\alpha}$.
Then, by Lemma 1.2,
there exists the induced automorphism $p$ of $\Lambda$
such that
$\vvp\big((F_{\bs q})_{\bs\alpha}\big)=Ft_{p(\bs\alpha)}$
for all $\bs\alpha\in\Lambda$.
Let
$\bs\sigma_i:=p^{-1}(\bs\varepsilon_{i})$
and
$x_i:=\vvp^{-1}(t_i)\in (F_{\bs q})_{\bs\sigma_i}$
for $i=1,\ldots,n$.
Then 
$\langle\bs\sigma_1,\ldots,\bs\sigma_n\rangle$
is a basis of $\Lambda$,
and
we have
$$x_{j}x_{i}
=\vvp^{-1}(t_j)\vvp^{-1}(t_i)
=\vvp^{-1}(t_jt_i)
=\vvp^{-1}(\eta_{ij}t_it_j)
=\eta_{ij}x_{i}x_{j}$$ 
for all $1\leq i<j\leq n$.
So (ii) holds.
Suppose that (ii) holds.
Since
$\langle\bs\sigma_{1},\ldots,\bs\sigma_{n}\rangle$
is a basis of $\Lambda$,
one has
$F_{\bs q}=\oplus_{\bs\alpha\in\Lambda}\ Fx_{\bs\alpha}$
where $x_{\bs\alpha}=x_1^{\alpha_1}\cdots x_n^{\alpha_n}$
for $\bs\alpha=\alpha_1\bs\sigma_1+\cdots +\alpha_n\bs\sigma_n$.
Define a map
$\vvp\:F_{\bs q}\longrightarrow
F_{\bs\eta}=F_{\bs\eta}[t_1^{\pm 1},\ldots,t_n^{\pm 1}]$
by $\vvp(x_{\bs\alpha})= t_{\bs\alpha}$
where $t_{\bs\alpha}=t_1^{\alpha_1}\cdots t_n^{\alpha_n}$
for all $\bs\alpha\in\Lambda$.
Then, since
$x_{j}x_{i}=\eta_{ij}x_{i}x_{j}$,
$\vvp$ is an isomorphism of algebras. 
Moreover, 
$\vvp$ is graded if we give
the 
$\langle\bs\sigma_{1},\ldots,\bs\sigma_{n}\rangle$-grading
to $F_{\bs\eta}$.
Hence (iii) holds.
Finally, (iii) clearly implies (i).
\qed
\enddemo

For convenience, we use the following notation:
\definition{Definition 1.4}
For square matrices $A_{1}$, $\ldots$, $A_{r}$
of sizes $l_{i}$, $i=1,\ldots,r$, we define the square matrix
$A_{1}\times\cdots\times A_{r}$ of size $l_{1}+\cdots+l_{r}$ to be
$$A_{1}\times\cdots\times A_{r}=
\left (\matrix
A_{1}  & \bk 1  & \bk 1       & \cdots &\bk 1 \\
\bk  1 & A_{2}  & \bk 1       &        &\vdots   \\
\bk 1  & \bk 1  & A_{3}       & \ddots &\vdots    \\
\vdots &        & \ddots      &\ddots  &\bk 1     \\
\bk 1  & \cdots & \cdots      & \bk 1  &  A_{r}          
\endmatrix
\right ),
$$
where $\bk 1$'s are matrices of suitable sizes whose entries are all $1$.
Also, we write $\bk 1=\bk 1_k$ if $\bk 1$ is a square matrix of size $k$.
\enddefinition 

\proclaim{Lemma 1.5}
\rom{(1)}
Let $\bs q=(q_{ij})$ be an $n\times n$ quantum matrix,
$\sigma$ a permutation on $\{1,\ldots,n\}$,
and put $\tilde{\bs q}_{\sigma}=(\tilde q_{ij})$
where $\tilde q_{ij}=q_{\sigma(i)\sigma(j)}$.
Then $\bs q\cong\tilde{\bs q}_{\sigma}$. 
In particular, for a transposition $(ij)\in S$,
we have $\bs q\cong\tilde{\bs q}_{(ij)}$. 

\rom{(2)} Let $\bs r$, $\bs s$
and $\bs\eta$ be quantum matrices
with $\bs s\cong \bs\eta$.
Then:
$$
\text{\rom{(i)}}\quad \bs r\times \bs s\cong\bs s\times \bs r,
\quad\quad\quad\quad\quad\quad
\text{\rom{(ii)}}\quad \bs r\times \bs s\cong\bs r\times \bs\eta.
$$ 
\endproclaim
\demo{Proof}
For (1),
let 
$F_{\bs q}=F_{\bs q}[t^{\pm 1}_{1},\ldots,t^{\pm 1}_{n}]$, and so we have $t_jt_i=q_{ij}t_it_j$.
Hence the generators
$\tilde t_i:=t_{\sigma(i)}$ satisfy
$\tilde t_j\tilde t_i=t_{\sigma(j)}t_{\sigma(i)}
=q_{\sigma(i)\sigma(j)}t_{\sigma(i)}t_{\sigma(j)}
=q_{\sigma(i)\sigma(j)}\tilde t_i\tilde t_j$,
and
$$F_{\bs q}=F_{\tilde{\bs q}_{\sigma}}
[\tilde t_1^{\pm 1},\ldots,
\tilde t_n^{\pm 1}].$$
Thus we get
$\bs q\cong\tilde{\bs q}_{\sigma}$.

For (2),
let $r$ and $s$ be the sizes of
the matrices $\bs r$ and $\bs s$,
respectively,
and 
let $n:=r+s$
and
$F_{\bs r\times \bs s}=
F_{\bs r\times \bs s}[t^{\pm 1}_{1},\ldots,t^{\pm 1}_{n}]$.

(i) follows from (1):
Take 
$$\sigma=\pmatrix
1&\cdots &s &s+1&\cdots &n \\
r+1 &\cdots&n&1&\cdots &r
\endpmatrix.
$$
Then
$\bs s\times \bs r=(\widetilde{\bs r\times \bs s})_{\sigma}$.

For (ii),
we consider a toral 
$\langle\bs\varepsilon_{1},\ldots,\bs\varepsilon_{n}\rangle$-grading 
of
$F_{\bs r\times \bs s}$.
Let $\bs r\times \bs\eta=(a_{ij})$.
The subalgebra of $F_{\bs r\times \bs s}$
generated by $t^{\pm 1}_{r+1},\ldots,t^{\pm 1}_{n}$
can be identified with
the quantum torus
$F_{\bs s}[t^{\pm 1}_{r+1},\ldots,t^{\pm 1}_{n}]$
with the
$\langle\bs\varepsilon_{r+1},\ldots,\bs\varepsilon_{n}\rangle$-grading.
By Lemma 1.3, our assumption $\bs s\cong\bs\eta$
implies that
there exists
a
basis 
$\langle\bs\sigma_{r+1},\ldots,\bs\sigma_{n}\rangle$
of $\Bbb Z\bs\varepsilon_{r+1}+\cdots +\Bbb Z\bs\varepsilon_{n}$
in $\Lambda$
such that
$x_{j}x_{i}=a_{ij}x_{i}x_{j}$
for all $r+1\leq i,j\leq n$
where $x_{i}$ is a nonzero element of degree $\bs\sigma_{i}$. 
Note that all $x_{1}:=t_{1},\ldots,x_{r}:=t_{r}$
commute with all $t_{r+1},\ldots,t_{n}$, and so all $x_{1},\ldots,x_{r}$
commute with all $x_{r+1},\ldots,x_{n}$.
Hence
we get 
$x_{j}x_{i}=a_{ij}x_{i}x_{j}$
for all $1\leq i,j\leq n$.
Since
$\langle\bs\varepsilon_{1},\ldots,\bs\varepsilon_{r},
\bs\sigma_{r+1},\ldots,\bs\sigma_{n}\rangle$
is a basis of $\Lambda$,
we obtain $\bs r\times \bs s\cong\bs r\times \bs\eta$
by Lemma 1.3.
\qed
\enddemo
\definition{Definition 1.6}
A quantum matrix $\bs\varepsilon=(\varepsilon_{ij})$
is called {\it elementary}
if
$\varepsilon_{ij}=1$ or $-1$ for all $i,j$.
Note that $\bs\varepsilon$ becomes a symmetric matrix.
Also,
the quantum torus $F_{\bs\varepsilon}$
determined by an elementary quantum matrix $\bs\varepsilon$
is called an {\it elementary quantum torus}.
\enddefinition
Note that any elementary quantum matrix is $\bk 1$ if $\ch F=2$.
Thus our argument will be trivial if $\ch F=2$, and so for convenience
we will assume that $\ch F\neq 2$ from now on.
\example{Example 1.7}
Let 
$$F_{\bs m_3}
=F_{\bs m_3}[t_1^{\pm 1},t_2^{\pm 1},t_3^{\pm 1}]
\quad\text{and}\quad
F_{\bs m_4}
=F_{\bs m_4}[t_1^{\pm 1},t_2^{\pm 1},t_3^{\pm 1},t_4^{\pm 1}]$$
be elementary quantum tori
with an
$\langle\bs\varepsilon_1,\bs\varepsilon_2,\bs\varepsilon_3\rangle$-grading
and
an $\langle\bs\varepsilon_1,\bs\varepsilon_2,
\bs\varepsilon_3,\bs\varepsilon_4\rangle$-grading,
respectively,
where
$$\bs m_3=
\left (\matrix \format\r&\quad\r &\quad\r    \\
   1   &-1    & -1   \\
  -1   & 1    & -1  \\
   -1   & -1   &   1   
\endmatrix
\right )
\quad\text{and}\quad
\bs m_4=
\left (\matrix \format\r&\quad\r &\quad\r &\quad\r   \\
   1   &-1    & -1   & -1 \\
  -1   & 1    & -1  & -1 \\
   -1   & -1   &   1  & -1 \\ 
  -1   & -1   &   -1  & 1
\endmatrix
\right ).
$$
In $F_{\bs m_3}$,
$t_1$ commutes with $t_2t_3$ which has degree 
$\bs\varepsilon_2+\bs\varepsilon_3$,
and
in $F_{\bs m_4}$,
$t_1$ commutes with $t_2t_3$ and $t_2t_4$ which has degree 
$\bs\varepsilon_2+\bs\varepsilon_3$
and
$\bs\varepsilon_2+\bs\varepsilon_4$.
Since 
$\langle\bs\varepsilon_1,\bs\varepsilon_2,
\bs\varepsilon_2+\bs\varepsilon_3\rangle$
and
$\langle\bs\varepsilon_1,\bs\varepsilon_2,
\bs\varepsilon_2+\bs\varepsilon_3,
\bs\varepsilon_2+\bs\varepsilon_4\rangle$
are bases of $\Lambda_3$ and $\Lambda_4$,
respectively,
we have by Lemma 1.3,
$$\bs m_3\cong
\left (\matrix \format\r&\quad\r &\quad\r    \\
   1   &-1    & 1   \\
  -1   & *    & *  \\
   1   & *    & * 
\endmatrix
\right )
\quad\text{and}\quad
\bs m_4\cong
\left (\matrix \format\r&\quad\r &\quad\r &\quad\r   \\
   1   &-1    & 1   & 1 \\
  -1   &  *   &  * & * \\
   1   & *  &   *  & * \\ 
  1    &  *  &  *  & *
\endmatrix
\right ),
$$
and the $*$-parts of both matrices are some elementary matrices.
Indeed in both algebras, we have $(t_2t_3)t_2=-t_2(t_2t_3)$,
and in $F_{\bs m_4}$, 
$(t_2t_4)t_2=-t_2(t_2t_4)$
and
$(t_2t_3)(t_2t_4)=-(t_2t_4)(t_2t_3)$.
So we get
$$\bs m_3\cong
\left (\matrix \format\r&\quad\r &\quad\r    \\
   1   &-1    & 1   \\
  -1   & 1    & -1  \\
   1   & -1    & 1 
\endmatrix
\right )
\quad\text{and}\quad
\bs m_4\cong
\left (\matrix \format\r&\quad\r &\quad\r &\quad\r   \\
   1   &-1    & 1   & 1 \\
  -1   &  1   &  -1 & -1 \\
   1   & -1  &   1  & -1 \\ 
  1    &  -1  &  -1  & 1
\endmatrix
\right ).
$$
In both algebras,
$t_1$ and $t_2$ commute with 
$t_1(t_2t_3)$ which has degree 
$\bs\varepsilon_1+\bs\varepsilon_2+\bs\varepsilon_3$,
and
in $F_{\bs m_4}$,
$t_1$ and $t_2$ commutes with 
$t_1(t_2t_4)$ which has degree 
$\bs\varepsilon_1+\bs\varepsilon_2+\bs\varepsilon_4$.
Since 
$\langle\bs\varepsilon_1,\bs\varepsilon_2,
\bs\varepsilon_1+\bs\varepsilon_2+\bs\varepsilon_3\rangle$
and
$\langle\bs\varepsilon_1,\bs\varepsilon_2,
\bs\varepsilon_1+\bs\varepsilon_2+\bs\varepsilon_3,
\bs\varepsilon_1+\bs\varepsilon_2+\bs\varepsilon_4\rangle$
are  bases of $\Lambda_3$ and $\Lambda_4$,
respectively,
we have by Lemma 1.3,
$$\bs m_3\cong
\left (\matrix \format\r&\quad\r &\quad\r    \\
   1   &-1    & 1   \\
  -1   & 1    & 1  \\
   1   & 1    & 1 
\endmatrix
\right )
\quad\text{and}\quad
\bs m_4\cong
\left (\matrix \format\r&\quad\r &\quad\r &\quad\r   \\
   1   &-1    & 1   & 1 \\
  -1   &  1   &  1 & 1 \\
   1   & 1  &   *  & * \\ 
  1    &  1  &  *  & *
\endmatrix
\right ),
$$
and
the $*$-part is $\bs h$
by
$(t_1t_2t_4)(t_1t_2t_3)
=-(t_1t_2t_3)(t_1t_2t_4)$.
Thus we have shown
$$
\bs m_3\cong\bs h_{1,3}
\quad\text{and}\quad 
\bs m_4\cong\bs h_{2,4}=\bs h\times\bs h.
$$
Note that we also have shown
$$
\align
F_{\bs m_3}
&\cong_{\Lambda}
F_{\bs h_{1,3}}[u_1^{\pm 1},u_2^{\pm 1},u_3^{\pm 1}]
\quad\text{via}\quad
t_1\mapsto u_1,\
t_2\mapsto u_2,\
t_1t_2t_3\mapsto u_3,\\
F_{\bs m_4}
&\cong_{\Lambda}
F_{\bs h\times\bs h}[u_1^{\pm 1},u_2^{\pm 1},
u_3^{\pm 1},u_4^{\pm 1}]
\quad\text{via}\quad
t_1\mapsto u_1,\
t_2\mapsto u_2,\
t_1t_2t_3\mapsto u_3\
t_1t_2t_4\mapsto u_4,
\endalign
$$
for the $\langle\bs\varepsilon_1,\bs\varepsilon_2,
\bs\varepsilon_3\rangle$-grading
of $F_{\bs h_{1,3}}$
and the
$\langle\bs\varepsilon_1,\bs\varepsilon_2,
\bs\varepsilon_3,\bs\varepsilon_4\rangle$-grading
of $F_{\bs h\times\bs h}$.
\endexample

\enskip

In general, the centre $Z(F_{\bs q})$ of a quantum torus
$F_{\bs q}$ is an algebra of Laurent polynomials, and the grading group
is given by
$$\{\bs\alpha\in\Lambda\ |\ \prod_{i,j}\ 
q_{ij}^{\alpha_{j}\beta_{i}}=1
\ \text{for all $\bs\beta\in\Lambda$}\}$$
(see [5] or [8]).
For later use,
we directly calculate the centre of $F_{\bs h_{l,n}}$.

\proclaim{Lemma 1.8}
Let $l>0$ and
$F_{\bs h_{l,n}}=F_{\bs h_{l,n}}[t^{\pm 1}_{1},\ldots,t^{\pm 1}_{n}]$
be an elementary torus.
Then the centre $Z(F_{\bs h_{l,n}})$ is equal to
$$F[t_{1}^{\pm 2},\ldots,t_{2l}^{\pm 2},t_{2l+1}^{\pm 1},\ldots,t_{n}^{\pm 1}],$$
the algebra of Laurent polynomials in the variables
$t_{1}^{2},\dots,t_{2l}^{2},t_{2l+1},\dots,
t_{n}$.
Hence for a 
$\langle\bs\sigma_1,
\ldots,\bs\sigma_n\rangle$-grading
of $F_{\bs h_{l,n}}$,
the grading group of $Z(F_{\bs h_{l,n}})$ is equal to
$$2\Bbb Z\bs\sigma_1+\cdots +2\Bbb Z\bs\sigma_{2l}
+\Bbb Z\bs\sigma_{2l+1}+\cdots +\Bbb Z\bs\sigma_n.$$
\endproclaim
\demo{Proof}
It is clear that
$Z':=F[t_{1}^{\pm 2},\ldots,t_{2l}^{\pm 2},t_{2l+1}^{\pm 1},\ldots,t_{n}^{\pm 1}]
\subset Z(F_{\bs h_{l,n}})=:Z$.
For the other inclusion,
if $Z\setminus Z'\neq\emptyset$,
there exists
$x:=t_{1}^{\kappa_1}\cdots t_{2l}^{\kappa_{2l}}\in Z$,
where $\kappa_i=0$ or $1$
but not all $\kappa_i$ are $0$.
But then, 
for $\kappa_j\neq 0$,
we have $xt_{k}=-t_{k}x$
where
$$k=
\cases
j+1&\text{if $j$ is odd} \\
j-1&\text{if $j$ is even}, 
\endcases
$$
i.e., $x\notin Z$,
which is a contradiction.
Hence $Z=Z'$.
\qed
\enddemo
Note that $\bs h_{0,n}=\bk 1$
and so 
$Z(F_{\bs h_{0,n}})=F[t^{\pm 1}_{1},\ldots,t^{\pm 1}_{n}]$.

\proclaim{Lemma 1.9}
Let $\bs\varepsilon=(\varepsilon_{ij})$
be an $n\times n$ elementary quantum matrix for $n\geq 3$.
If $\varepsilon_{kp}=\varepsilon_{kq}=-1$
for some distinct $1\leq k,p,q\leq n$,
then there exists
an elementary quantum matrix
$\bs\eta=(\eta_{ij})$ with 
$$\align
\eta_{ij}&=\varepsilon_{ij}\quad\text{for all $i,j\neq q$}
\quad (\eta_{qq}=\varepsilon_{qq}=1), \\
\eta_{iq}&=\varepsilon_{ip}\varepsilon_{iq}
\quad\text{for all $i\neq q$}
\endalign
$$
such that
$\bs\varepsilon\cong\bs\eta$.
In particular,

\rom{(a)}
$\eta_{kq}=1$
and
$\eta_{ki}=\varepsilon_{ki}$ for all $i\neq q$.

\rom{(b)}
if 
$k=2$
and
$p=1$,
then 
$\eta_{i1}=\varepsilon_{i1}$ for all $i$,
i.e., the first rows of $\bs\varepsilon$ and $\bs\eta$
are the same.
\endproclaim
\demo{Proof}
Let 
$F_{\bs\varepsilon}=F_{\bs\varepsilon}[t^{\pm 1}_{1},\ldots,t^{\pm 1}_{n}]$
with a
$\langle\bs\sigma_{1},\ldots,\bs\sigma_{n}\rangle$-grading.
Since $\varepsilon_{kp}=\varepsilon_{kq}=-1$,
we have
$t_pt_k=-t_kt_p$ and
$t_qt_k=-t_kt_q$.
Hence $t_k$ commutes with $t_pt_q$
which has degree $\bs\sigma_p+\bs\sigma_q$.
Let
$$x_{1}:=t_{1},\ \ldots,
x_{q-1}:=t_{q-1},\
x_{q}:=t_{p}t_{q},\
x_{q+1}:=t_{q+1},\ \ldots,\ x_{n}:=t_{n}.$$
Then the relations between $x_i$ and $x_j$
for $1\leq i,j\leq n$
determine an elementary quantum matrix
$\bs\eta=(\eta_{ij})$,
i.e.,
$x_jx_i=\eta_{ij}x_ix_j$.
It is clear that
$\eta_{ij}=\varepsilon_{ij}$
for all $i,j\neq q$.
For $i\neq q$, we have
$x_qx_i=(t_pt_q)t_i=\varepsilon_{ip}\varepsilon_{iq}t_i(t_pt_q)
=\varepsilon_{ip}\varepsilon_{iq}x_ix_q$.
Hence $\eta_{iq}=\varepsilon_{ip}\varepsilon_{iq}$.
Since
$$\langle\bs\sigma_{1},\ldots,
\bs\sigma_{q-1},\bs\sigma_{p}+\bs\sigma_{q},
\bs\sigma_{q+1},
\ldots,\bs\sigma_{n}\rangle$$ 
is a basis of $\Lambda$,
we get $\bs\varepsilon\cong\bs\eta$ by Lemma 1.3.
(a) and (b) are clear now.
\qed
\enddemo
Our first result is the following:
\proclaim{Theorem 1.10}
Let $\bs\varepsilon$ 
be an $n\times n$ elementary quantum matrix.
Then 
there exists $l\geq 0$ such that
 $\bs\varepsilon\cong\bs h_{l,n}$ where
$$\bs h_{l,n}=
\overbrace
{\bs h\times\cdots\times\bs h}^{\text{$l$-times}}
\times 
\bk 1_{n-2l}
\quad\text{and}\quad
\bs h=
\left (\matrix \format\r&\quad\r \\
    1   & -1      \\
   -1   &  1          
\endmatrix
\right ).
$$
Also, there exists a 
$\langle\bs\sigma_1,
\ldots,\bs\sigma_n\rangle$-grading
of $F_{\bs\varepsilon}$ such that
the grading group of the centre $Z(F_{\bs\varepsilon})$ is equal to
$$2\Bbb Z\bs\sigma_1+\cdots +2\Bbb Z\bs\sigma_{2l}
+\Bbb Z\bs\sigma_{2l+1}+\cdots +\Bbb Z\bs\sigma_n.$$
Moreover, the number $l$ is an isomorphism invariant of $F_{\bs\varepsilon}$.
\endproclaim
\demo{Proof}
We prove this by induction on $n$.
When $n=1$, $\bs\varepsilon$ has to be $(1)$, and so the statement is clear.
Let $n>1$,
$\bs\varepsilon=(\varepsilon_{ij})$
and 
$$N_k(\bs\varepsilon):
=|\{i\ |\ \varepsilon_{ki}=-1,\ 1\leq i\leq n\}|$$
where $|\ |$ is the number of elements of a set.
(We will use this notation only for $k=1$ and $2$.)
If $N_1(\bs\varepsilon)=0$,
then $\bs\varepsilon=(1)\times \bs\varepsilon'$
for an
elementary quantum matrix
$\bs\varepsilon'$
of size $n-1$.
By induction,
we have
$\bs\varepsilon'\cong\bs h_{l,n-1}$
for some $l\geq 0$.
Then, by Lemma 1.5(2),
we get 
$$\bs\varepsilon=(1)\times \bs\varepsilon'\cong
(1)\times \bs h_{l,n-1}\cong
\bs h_{l,n-1}\times (1)
=\bs h_{l,n}.$$

If $N_1(\bs\varepsilon)>1$,
then by Lemma 1.9(a) for $k=1$,
there exists an elementary quantum matrix
$\bs\varepsilon'$ such that
$\bs\varepsilon\cong\bs\varepsilon'$ and
$N_1(\bs\varepsilon')=N_1(\bs\varepsilon)-1$.
Repeating this, 
we obtain an elementary quantum matrix 
$\bs\nu$
such that 
$\bs\varepsilon\cong\bs\nu$
and $N_1(\bs\nu)=1$,
i.e.,
only one entry, say
the $(1i_0)$-entry, is $-1$ in the first column of $\bs\nu$.
So if $N_1(\bs\varepsilon)=1$,
we can also put $\bs\nu=\bs\varepsilon$.
Then, by Lemma 1.5(1),
we get
$$\bs\varepsilon\cong
\bs\nu_{(2i_0)}=:
\bs\eta=(\eta_{ij})=
\left (\matrix
   1   & -1  & 1 &\cdots & 1\\
  -1   &     &   &       &  \\
   1   &     & * &       &  \\
\vdots &     &   &       &  \\
   1   &     &   &       &
\endmatrix
\right ),
$$
$\eta_{12}=\eta_{21}=-1$,
the other $\eta_{1i}=\eta_{i1}=1$
and $*$ is some elementary quantum matrix
of size $n-1$.

If $n=2$, we have
$\bs\eta=\bs h$ and we are done.
We assume that  $n>2$.
Note that $N_2(\bs\eta)\geq 1$
since $\eta_{21}=-1$.
If $N_2(\bs\eta)> 1$,
we can apply
Lemma 1.9(b) for any $q>2$ such that $\eta_{2q}=-1$,
and get an elementary quantum matrix 
$\bs\eta'$ such that 
$\bs\eta\cong \bs\eta'$,
$N_1(\bs\eta')=N_1(\bs\eta)=1$ and
$N_2(\bs\eta')=N_2(\bs\eta)-1$.
Repeating this,
we obtain an elementary quantum matrix 
$\bs\mu=(\mu_{ij})$ such that 
$\bs\eta\cong \bs\mu$,
$N_1(\bs\mu)=N_2(\bs\mu)=1$
and
$\mu_{21}=\mu_{12}=-1$.
Also, if $N_2(\bs\eta)=1$,
we put $\bs\eta=\bs\mu$.
Thus we have
$\bs\eta\cong \bs\mu=
\bs h\times\bs\mu'$
for an elementary quantum matrix $\bs\mu'$
of size $n-2$.
By induction,
we have
$\bs\mu'\cong\bs h_{l',n-2}$
for some $l'\geq 0$.
Then, by Lemma 1.5(2)(ii),
we get 
$\bs\mu=\bs h\times \bs\mu'\cong
\bs h\times \bs h_{l',n-2}
=\bs h_{l,n}$
where $l=l'+1$, and hence 
$\bs\varepsilon\cong\bs\eta
\cong\bs\mu
\cong\bs h_{l,n}$.

The description of the centre follows from Lemma 1.3 and
Lemma 1.8.
For the last statement,
suppose that $\bs h_{l,n}\cong \bs h_{l',n}$.
Then, by Lemma 1.3,
$F_{\bs h_{l,n}}\cong_{\Lambda}F_{\bs h_{l',n}}$ for some toral gradings.
Hence the grading groups of the centres 
of $F_{\bs h_{l,n}}$ and $F_{\bs h_{l',n}}$
coincide,
which implies $l=l'$ by Lemma 1.8.
Therefore, $l$ is an isomorphism invariant of $F_{\bs\varepsilon}$.
\qed
\enddemo

\head
\S 2  Elementary quantum tori with graded involution
\endhead

From now on, we always consider a quantum torus
as a toral $\Lambda$-graded algebra.
Let $F_{\bs q}=F_{\bs q}
[t^{\pm 1}_{1},
\ldots,t_{n}^{\pm 1}]$ 
be the quantum torus determined by $\bs q=(q_{ij})$,
and let $\tau$ be a graded involution of 
$F_{\bs q}$.
Then we have $\tau(t_{i})=a_{i}t_{i}$ for 
some $a_{i}\in F$, $i=1,\ldots,n$.
Since $t_{i}=\tau^2(t_{i})=a_{i}^2t_{i}$,
one gets $a_{i}=\pm 1$ for all $1\leq i\leq n$.
Moreover, one has
$$
a_{i}a_{j}q_{ij}t_{j}t_{i}
=\tau(q_{ij}t_{i}t_{j})
=\tau(t_{j}t_{i})=a_{i}a_{j}t_{i}t_{j}
=a_{i}a_{j}q_{ji}t_{j}t_{i},$$
and hence $q_{ij}^{-1}=q_{ji}$,
i.e.,
$q_{ij}=\pm 1$
for all $1\leq i,j\leq n$.
Thus 
$\bs q$ has to be elementary.

Conversely, 
it is straightforward to check that
for an elementary quantum tours
$F_{\bs\varepsilon}=F_{\bs\varepsilon}
[t^{\pm 1}_{1},
\ldots,t_{n}^{\pm 1}]$  and
each $(a_{1},\ldots,a_{n})$,
$a_{i}=\pm 1$,
there exists a unique involution of
$F_{\bs\varepsilon}$ such that
$\tau(t_{i})=a_{i}t_{i}$ for all $1\leq i\leq n$.
We call this $\tau$
{\it of type $(a_{1},\ldots,a_{n})$},
denoted
$\tau=(a_{1},\ldots,a_{n})$.
The graded involution of type
$(1,\ldots,1)$ is called
the {\it main involution},
denoted $*$.
Thus we have the following proposition, which is stated in [1]:
\proclaim{Proposition 2.1}
Let $F_{\bs q}=F_{\bs q}[t^{\pm 1}_{1},\ldots,t^{\pm 1}_{n}]$
be a quantum torus over $F$.
Then there exists a graded involution $\tau$ of
$F_{\bs q}$
if and only if
$\bs q$ is elementary.
In this case, 
$\tau$ has type $(a_{1},\ldots,a_{n})$,
i.e., $\tau(t_{i})=a_{i}t_{i}$
where $a_{i}=1$ or $-1$
for all $1\leq i\leq n$.
\qed
\endproclaim
Recall the notion of isomorphism in the class of
algebras with involution.
Namely, for algebras with involution $(A,\tau)$ and $(B,\rho)$,
 an {\it isomorphism
of algebras with involution} from $(A,\tau)$ onto $(B,\rho)$
is
an isomorphism $f$ from $A$ onto $B$
satisfying $f\tau=\rho f$,
and in this case
we denote this by $(A,\tau)\cong(B,\rho)$.
Moreover, if $A$ and $B$ are $\Lambda$-graded algebras,
$\tau$ and $\rho$ are graded involutions
and the $f$ happens to be a graded isomorphism,
we write $(A,\tau)\cong_{\Lambda}(B,\rho)$.
Finally, the {\it centre} $Z(A,\tau)$ of $(A,\tau)$ is defined as
$$Z(A,\tau)=Z(A)\cap \{a\in A\ |\ \tau(a)=a\},$$
where $Z(A)$ is the centre of the algebra $A$.

One can prove the following lemmas similar to Lemma 1.3 and 1.5.
Since the proofs can be done in the same manner,
they will be left to the reader.
\proclaim{Lemma 2.2}
Let
$(F_{\bs\varepsilon},\tau)$ 
and $(F_{\bs\eta},\rho)$ be
elementary quantum tori with graded involutions.
Let $\bs\eta=(\eta_{ij})_{1\leq i,j\leq n}$ and $\rho=(a_{1},\ldots,a_{n})$.
Then the following are equivalent:

\rom{(i)} $(F_{\bs\varepsilon},\tau)\cong (F_{\bs\eta},\rho)$, 

\rom{(ii)} for any toral grading of $F_{\bs\varepsilon}$,
there exists a basis
$\langle\bs\sigma_{1},\ldots,\bs\sigma_{n}\rangle$ 
of $\Lambda$ and
nonzero homogeneous elements $x_{i}
\in F_{\bs\varepsilon}$ of degree $\bs\sigma_{i}$
such that
$x_{j}x_{i}=\eta_{ij}x_{i}x_{j}$ and $\tau(x_{i})=a_{i}x_{i}$
for all $1\leq i<j\leq n$,

\rom{(iii)} for any toral grading of $F_{\bs\varepsilon}$,
there exists a toral grading
of $F_{\bs\eta}$
such that
$(F_{\bs\varepsilon},\tau)\cong_{\Lambda} (F_{\bs\eta},\rho)$. 
In that case,
the grading subgroups of the centres $Z(F_{\bs\varepsilon},\tau)$ and 
$Z(F_{\bs\eta},\rho)$
coincide.
\qed
\endproclaim
For graded involutions $\tau$ and $\rho$
of type $(a_{1},\ldots, a_r)$ and $(b_{1},\ldots, b_s)$, respectively,
we denote the graded involution
of type $(a_{1},\ldots, a_r,b_{1},\ldots, b_s)$
by $\tau\times \rho$.

\proclaim{Lemma 2.3}
Let $(F_{\bs r},\tau)$, $(F_{\bs s},\rho)$
and $(F_{\bs\eta},\rho_1)$ be elementary quantum tori
with graded involutions.
Assume that $(F_{\bs s},\rho)\cong (F_{\bs\eta},\rho_1)$.
Then:
$$
\text{\rom{(i)}}\quad (F_{\bs r\times \bs s},\tau\times\rho)\cong 
(F_{\bs s\times \bs r},\rho\times\tau),
\quad\quad
\text{\rom{(ii)}}\quad (F_{\bs r\times \bs s},\tau\times\rho)\cong
(F_{\bs r\times \bs\eta},\tau\times\rho_1).
\qed
$$ 
\endproclaim

\enskip

We start to classify elementary tori with graded involution.
Let $\tau$ be a graded involution of an elementary quantum torus
$F_{\bs\varepsilon}$.
Then, by Theorem 1.10 and Lemma 1.3,
we have $F_{\bs\varepsilon}\cong_{\Lambda} F_{\bs h_{l,n}}$
for some $l\geq 0$
and toral gradings,
and hence $(F_{\bs\varepsilon},\tau)\cong_{\Lambda} (F_{\bs h_{l,n}},\rho)$
for some graded involution $\rho$ of $F_{\bs h_{l,n}}$.
Thus it is enough to classify $F_{\bs h_{l,n}}$
with graded involutions.
Besides the main involution 
$*=(1,\ldots,1)$, 
we define two specific graded involutions
of $F_{\bs h_{l,n}}$, namely,
$$
\align
\tau_1
&=(1,\ldots,1,-1,1,\ldots,1), \\
&\text{where only the $2l+1$ position is $-1$,
if $n-2l\geq 1$} \\
\tau_2
&=(1,\ldots,1,-1,-1,1,\ldots,1), \\
&\text{where only the $2l-1$ and $2l$ positions are $-1$,
if $l\geq 1$.}
\endalign
$$

\remark{Remark}
By Lemma 1.8, $*$ and $\tau_2$ fix the centre $Z$ of $F_{\bs h_{l,n}}$
but $\tau_1$ does not.  It is easily seen that the {\it central closure}
$\ol F_{\bs h_{l,n}}=\ol Z\otimes_Z F_{\bs h_{l,n}}$
is a simple algebra over $\ol Z$, where $\ol Z$ is the field of fractions
of $Z$.
Let $\tau=*$, $\tau_1$ or $\tau_2$.
By the universal property of
the central closure $\ol F_{\bs h_l^{(n)}}$,
the natural extension $\ol{\tau}$ of
$\tau$ 
 defined by
$\ol{\tau}(z\otimes x)=\tau(z)\otimes \tau(x)$
is an involution of $\ol F_{\bs h_l^{(n)}}$.
Since $\ol *$ and $\ol{\tau_2}$
fix $\ol Z$,
they are involutions of first kind,
while $\ol{\tau_1}$ does not, and so it is an involution of second kind.
\endremark

\example{Example 2.4}
Recall the two elementary quantum matrices
$\bs m_3$ and $\bs m_4$
defined in Example 1.7.
The isomorphisms $\bs m_3\cong \bs h_{1,3}$
and $\bs m_4\cong \bs h_{2,4}$ there give isomorphisms of algebras
with involution, namely,
$$(F_{\bs m_3},*)\cong (F_{\bs h_{1,3}},\tau_1)
\quad\text{and}\quad
(F_{\bs m_4},*)\cong (F_{\bs h_{2,4}},\tau_2).$$
\endexample
Like Lemma 1.8, we have the following lemma about the centres:

\proclaim{Lemma 2.5}
Let 
$F_{\bs h_{l,n}}=F_{\bs h_{l,n}}[t^{\pm 1}_{1},\ldots,t^{\pm 1}_{n}]$
be an elementary torus.
Then 
$$\align
Z(F_{\bs h_{l,n}},*)=Z(F_{\bs h_{l,n}},\tau_2)
&=F[t_{1}^{\pm 2},\ldots,t_{2l}^{\pm 2},t_{2l+1}^{\pm 1},\ldots,t_{n}^{\pm 1}] \\
Z(F_{\bs h_{l,n}},\tau_1)
&=F[t_{1}^{\pm 2},\ldots,t_{2l+1}^{\pm 2},t_{2l+2}^{\pm 1},\ldots,t_{n}^{\pm 1}].
\endalign
$$
(For $(F_{\bs h_{l,n}},\tau_2)$, we are always assuming $l\geq 1$, but for the others,
$l$ can be $0$.)

Hence for a 
$\langle\bs\sigma_1,
\ldots,\bs\sigma_n\rangle$-grading
of $F_{\bs h_{l,n}}$,
the grading groups of $Z(F_{\bs h_{l,n}},*)$ 
and $Z(F_{\bs h_{l,n}},\tau_2)$ are equal to
$$2\Bbb Z\bs\sigma_1+\cdots +2\Bbb Z\bs\sigma_{2l}
+\Bbb Z\bs\sigma_{2l+1}+\cdots +\Bbb Z\bs\sigma_n,$$
and the grading group of $Z(F_{\bs h_{l,n}},\tau_1)$ is equal to
$$2\Bbb Z\bs\sigma_1+\cdots +2\Bbb Z\bs\sigma_{2l+1}
+\Bbb Z\bs\sigma_{2l+2}+\cdots +\Bbb Z\bs\sigma_n.$$
\endproclaim
\demo{Proof}
From Lemma 1.8, we already knows the description of the centre
$Z(F_{\bs h_{l,n}})$
of $F_{\bs h_{l,n}}$.
So only the fixed elements of $Z(F_{\bs h_{l,n}})$
under each $*$, $\tau_1$ and $\tau_2$
have to be calculated.
This easy exercise is left to the reader.
\qed
\enddemo

For the classification of elementary tori with graded involution,
we use the following:
\proclaim{Lemma 2.6}
Let
$*$ be the main involution
and $\tau_1$ the graded involution
of $F_{\bs h_{l,n}}$ defined above.
Then:
\roster
\item"(i)"
$\big(F_{\bs h},(1,-1)\big)\cong 
\big(F_{\bs h},(-1,1)\big)\cong(F_{\bs h},*)$,
\item"(ii)"
$\big(F_{\bk 1_2},(-1,-1)\big)\cong(F_{\bk 1_2},\tau_1)$,
\item"(iii)"
$\big(F_{\bs h_{1,3}},(-1,-1,-1)\big) 
\cong(F_{\bs h_{1,3}},\tau_1)$,
\item"(iv)"
$\big(F_{\bs h_{2,4}},(-1,-1,-1,-1)\big) 
\cong(F_{\bs h_{2,4}},*)$.
\endroster
\endproclaim
\demo{Proof}
Let 
$F_{\bs h_{l,n}}=F_{\bs h_{l,n}}[t^{\pm 1}_{1},\ldots,t^{\pm 1}_{n}]$
with an
$\langle\bs\varepsilon_{1},\ldots,\bs\varepsilon_{n}\rangle$-grading.
Then we note that
$t_{i_1}\cdots t_{i_r}$ has degree 
$\bs\varepsilon_{i_1}+\cdots +\bs\varepsilon_{i_r}$.

For (i), we have 
$n=2$ and $l=1$.
Let $\tau=(1,-1)$.
Then we have $\tau(t_{1})=t_{1}$ and $\tau(t_{2})=-t_{2}$.
Since
$(t_{1}t_{2})t_{1}=-t_{1}(t_{1}t_{2})$
and $\tau(t_{1}t_{2})=t_{1}t_{2}$,
and since $\langle\bs\varepsilon_{1},
\bs\varepsilon_{1}+\bs\varepsilon_{2}\rangle$
is a basis of $\Lambda_{2}$,
we get $(F_{\bs h},\tau)\cong (F_{\bs h},*)$
by Lemma 2.2.
The case $(-1,1)$ can be proven in the same way.

For (ii), we have 
$n=2$ and $l=0$.
Let $\tau=(-1,-1)$. 
Then we have
$\tau(t_{1})=-t_{1}$ and $\tau(t_{2})=-t_{2}$.
Since
$(t_{1}t_{2})=t_{1}(t_{1}t_{2})$
and
$\tau(t_{1}t_{2})=t_{1}t_{2}$, 
and since
$\langle\bs\varepsilon_{1},
\bs\varepsilon_{1}+\bs\varepsilon_{2}\rangle$
is a basis of $\Lambda_{2}$,
we get $(F_{\bk 1_2},\tau)\cong (F_{\bk 1_2},\tau_1)$
by Lemma 2.2.

For (iii), we have 
$n=3$ and $l=1$.
Let $\tau=(-1,-1,-1)$. 
Then we have
$\tau(t_{1})=-t_{1}$, $\tau(t_{2})=-t_{2}$ and $\tau(t_{3})=-t_{3}$.
Since
$(t_{2}t_{3})(t_{1}t_{2}t_{3})=-(t_{1}t_{2}t_{3})(t_{2}t_{3})$,
$t_{3}(t_{1}t_{2}t_{3})=(t_{1}t_{2}t_{3})t_{3}$,
$t_{3}(t_{2}t_{3})=(t_{2}t_{3})t_{3}$,
$\tau(t_{1}t_{2}t_{3})=t_{1}t_{2}t_{3}$
and $\tau(t_{2}t_{3})=t_{2}t_{3}$,
and since $\langle\bs\varepsilon_{1}+\bs\varepsilon_{2}+\bs\varepsilon_{3},
\bs\varepsilon_{2}+\bs\varepsilon_{3},
\bs\varepsilon_{3}\rangle$
is a basis of $\Lambda_{3}$,
we get $(F_{\bs h_{1,3}},\tau)\cong (F_{\bs h_{1,3}},\tau_1)$
by Lemma 2.2.

For (iv), we have 
$n=4$ and $l=2$.
Let $\tau=(-1,-1,-1,-1)$. 
Then we have
$\tau(t_{1})=-t_{1}$, $\tau(t_{2})=-t_{2}$,
$\tau(t_{3})=-t_{3}$ and $\tau(t_{4})=-t_{4}$.
Put
$x_{1}:=t_{1}t_{2}t_{4}$,
$x_{2}:=t_{2}t_{4}$,
$x_{3}:=t_{1}t_{3}$
and
$x_{4}:=t_{1}t_{3}t_{4}$.
Then one can check that
$x_{j}x_{i}=a_{ij}x_{i}x_{j}$
where $(a_{ij})=\bs h_{2,4}$
and $\tau(x_{i})=x_{i}$
for $1\leq i,j\leq 4$.
Also, one can check that
$$\langle\bs\varepsilon_{1}+
\bs\varepsilon_{2}+\bs\varepsilon_{4},
\bs\varepsilon_{2}+\bs\varepsilon_{4},
\bs\varepsilon_{1}+\bs\varepsilon_{3},
\bs\varepsilon_{1}+
\bs\varepsilon_{3}+\bs\varepsilon_{4}\rangle$$
is a basis of $\Lambda_{4}$.
Hence by Lemma 2.2,
we get $(F_{\bs h_{2,4}},\tau)
\cong (F_{\bs h_{2,4}},*)$.
\qed
\enddemo
Now we state one of our main theorems.
\proclaim{Theorem 2.7}
Let $\tau$
be an arbitrary graded involution of 
an elementary quantum torus $F_{\bs\varepsilon}$.
Let
$*$ be the main involution,
and $\tau_1$ and $\tau_2$ the graded involutions
of $F_{\bs h_{l,n}}$ defined above.
Then
$(F_{\bs\varepsilon},\tau)$ is graded isomorphic to exactly one of
$$(F_{\bs h_{l,n}},*), \quad
(F_{\bs h_{l,n}},\tau_{1})\quad\text{or}\quad
(F_{\bs h_{l,n}},\tau_{2}),$$
and for each of these $l$ is an invariant of the isomorphism class.
Moreover, we have
\roster
\item"(i)"
$(F_{\bs\varepsilon},*)\cong (F_{\bs h_{l,n}},\tau_{1})
\Longrightarrow l\geq 1$;
\item"(ii)"
$(F_{\bs\varepsilon},*)\cong (F_{\bs h_{l,n}},\tau_{2})
\Longrightarrow l\geq 2$;
\item"(iii)"
$(F_{\bs h_{l,n}},\tau_{1})
\cong_{\Lambda} (F_{\bs h_{l-1,n-3}\times \bs m_3},*)$
for $l\geq 1$;
\item"(iv)"
$(F_{\bs h_{l,n}},\tau_{2})
\cong_{\Lambda} (F_{\bs h_{l-2,n-4}\times \bs m_4},*)$
for $l\geq 2$,
where $\bs m_3$ and $\bs m_4$
are the elementary quantum matrices defined in Example 1.7.
\endroster
In particular, $(F_{\bs\varepsilon},\tau)$ 
is graded isomorphic to exactly one of 
$(F_{\bs h_{0,n}},\tau_{1})$,
$(F_{\bs h_{1,n}},\tau_{2})$ or
$(F_{\bs\eta},*)$
for some elementary quantum matrix $\bs\eta$.
\endproclaim
\demo{Proof}
We have $(F_{\bs\varepsilon},\tau)\cong_{\Lambda} (F_{\bs h_{l,n}},\rho)$
for some graded involution $\rho$ of $F_{\bs h_{l,n}}$
as mentioned above.
So we classify $(F_{\bs h_{l,n}},\rho)$ for
$\rho=(a_{1},\ldots,a_{n})$.
Note that $\bs h_{l,n}=\bs h_{l,2l}\times \bk 1_{n-2l}$.
We consider
$\big(F_{\bs h_{l,2l}},(a_{1},\ldots,a_{2l})\big)$
and
$\big(F_{\bk 1_{n-2l}},(a_{2l+1},\ldots,a_{n})\big)$
separately.
By Lemma 2.3
and Lemma 2.6(i) and (iv),
we have
$$\big(F_{\bs h_{l,2l}},(a_{1},\ldots,a_{2l})\big)
\cong(F_{\bs h_{l,2l}},*)
\quad\text{or}\quad
(F_{\bs h_{l,2l}},\tau_2),$$
and by Lemma 2.6(ii),
$$\big(F_{\bk 1_{n-2l}},(a_{2l+1},\ldots,a_{n})\big)
\cong (F_{\bk 1_{n-2l}},*)
\quad\text{or}\quad
(F_{\bk 1_{n-2l}},\tau_1).$$
Hence by Lemma 2.3, we get
$$(F_{\bs h_{l,n}},\rho)
\cong (F_{\bs h_{l,n}},*),\
(F_{\bs h_{l,n}},\tau_{1}),\
(F_{\bs h_{l,n}},\tau_{2})
\ \text{or}\ 
\big(F_{\bs h_{l,n}},(1,\ldots,1,-1,-1,-1,1,\ldots,1)\big),$$
and the last one is isomorphic to 
$(F_{\bs h_{l,n}},\tau_{1})$
by Lemma 2.6(iii).
Hence, by
Lemma 2.2,
we have obtained
$(F_{\bs\varepsilon},\tau)\cong_{\Lambda} (F_{\bs h_{l,n}},*)$,
$(F_{\bs h_{l,n}},\tau_{1})$
or $(F_{\bs h_{l,n}},\tau_{2})$.

By Lemma 2.5, we know
the grading groups of the centres $Z(F_{\bs h_{l,n}},*)$,
$Z(F_{\bs h_{l,n}},\tau_{1})$ and $Z(F_{\bs h_{l,n}},\tau_{2})$,
and hence by Lemma 2.2,
$l$ is an invariant of the isomorphism classes.
Moreover, 
the grading groups of $Z(F_{\bs h_{l,n}},*)$
and $Z(F_{\bs h_{l,n}},\tau_{2})$,
are the same but different from the one of $Z(F_{\bs h_{l,n}},\tau_{1})$.
Thus, by Lemma 2.2, we get
$(F_{\bs h_{l,n}},*)\not\cong (F_{\bs h_{l,n}},\tau_{1})$
and
$(F_{\bs h_{l,n}},\tau_{1})\not\cong (F_{\bs h_{l,n}},\tau_{2})$.
We postpone the proof of
$(F_{\bs h_{l,n}},*)\not\cong (F_{\bs h_{l,n}},\tau_{2})$
until \S 4
(right after the proof of Lemma 4.1).

(i):
Suppose that $(F_{\bs h_{0,n}},\tau_{1})\cong (F_{\bs\varepsilon},*)$.
We have
$\bs h_{0,n}=\bk 1$, which forces $\bs\varepsilon=\bk 1$,
and hence $*$ is the identity map.
This is a contradiction
since $\tau_1$ is not
the identity map.
Therefore, we get $(F_{\bs h_{0,n}},\tau_{1})\ncong (F_{\bs\varepsilon},*)$.

(ii):
Suppose that $(F_{\bs h_{1,n}},\tau_{2})\cong (F_{\bs\varepsilon},*)$.
Let
$F_{\bs h_{1,n}}=F_{\bs h_{1,n}}[t^{\pm 1}_{1},\ldots,t^{\pm 1}_{n}]$ 
with an $\langle\bs\varepsilon_{1},\ldots,\bs\varepsilon_{n}\rangle$-grading.
By Lemma 2.2,
there exists a basis
$\langle\bs\rho_{1},\ldots,\bs\rho_{n} \rangle$ of $\Lambda$
such that
a nonzero element $x_i\in F_{\bs h_{1,n}}$ of degree $\bs\rho_i$ 
are fixed by $\tau_{2}$
for all $i=1,\ldots,n$. 
Let
$\bs\rho_{i}=\alpha_{i1}\bs\varepsilon_{1}+\cdots+\alpha_{in}\bs\varepsilon_{n}$
for
$\alpha_{ij}\in\Bbb Z$.
Then one can take 
$x_i=t^{\alpha_{i1}}_{1}\cdots t_{n}^{\alpha_{in}}$.
Since $\tau_{2}=(-1,-1,\ 1,\ldots,1)$,
we have,
by the multiplication rule (1.1) of a quantum torus,
$$\tau_{2}(x_i)
=(-1)^{\alpha_{i1}+\alpha_{i2}}
t^{\alpha_{in}}_{1}\cdots t_{n}^{\alpha_{i1}}
=(-1)^{\alpha_{i1}+\alpha_{i2}+\alpha_{i1}\alpha_{i2}}
x_i
=x_i.$$
Hence $\alpha_{i1}$ and $\alpha_{i2}$ are both even for
all $i=1,\ldots,n$.
This implies that
the determinant of the matrix $(\alpha_{ij})$
is even.
This is absurd
since
$\langle\bs\rho_{1},\ldots,\bs\rho_{n} \rangle$ 
is a basis of $\Lambda$.
Therefore, we get $(F_{\bs h_{1,n}},\tau_{2})\ncong (F_{\bs\varepsilon},*)$.

For (iii) and (iv),
let
$F_{\bs h_{l,n}}=F_{\bs h_{l,n}}[t^{\pm 1}_{1},\ldots,t^{\pm 1}_{n}]$.
Let $U$ be the subalgebra of $(F_{\bs h_{l,n}},\tau_1)$
generated by $t_{2l-1}^{\pm 1}$, $t_{2l}^{\pm 1}$ and $t_{2l+1}^{\pm 1}$,
and let
$V$ be the subalgebra of $(F_{\bs h_{l,n}},\tau_2)$
generated by $t_{2l-3}^{\pm 1}$, $t_{2l-2}^{\pm 1}$
$t_{2l-1}^{\pm 1}$ and $t_{2l}^{\pm 1}$.
Then we have 
$(U,\tau_1\mid_U)\cong (F_{\bs h_{1,3}},\tau_1)\cong (F_{\bs m_3},*)$
and
$(U,\tau_2\mid_V)\cong (F_{\bs h_{2,4}},\tau_2)\cong (F_{\bs m_4},*)$
(see Example 2.4).
Therefore, by Lemma 2.3
we obtain
(iii) and (iv).
\qed
\enddemo

\head
\S 3  Semilattices
\endhead

We review semilattices (see [1]).
Let $\Bbb E$ be a Euclidean space.
A subset $S$ of $\Bbb E$ is called a {\it semilattice in $\Bbb E$} if
\roster
\item"(S1)" $0\in S$,
\item"(S2)" $S-2S\subset S$,
\item"(S3)" $S$ spans $\Bbb E$,
\item"(S4)" $S$ is discrete in $\Bbb E$.
\endroster
Also, a subset $S$ of a free abelian group
of finite rank is called a {\it semilattice in $\Lambda$} if
(S1), (S2) and
\roster
\item"(S3)$'$" $S$ spans $\Lambda$.
\endroster
If $S$ is a semilattice in $\Bbb E$, then
the group $\langle S\rangle$ generated by $S$
is a lattice in $\Bbb E$ and
$S$ is a semilattice in $\langle S\rangle$.
Also, if $S$ is a semilattice in $\Lambda$, then
$S$ can be considered  as a semilattice in some $\Bbb E$.
Note that $2S$ is not a semilattice in $\langle S\rangle$,
but a semilattice in $\Bbb E$.
We define the {\it rank of a semilattice $S$ in $\Bbb E$}
(resp. {\it in $\Lambda$}) as the dimension of $\Bbb E$
(resp. the rank of $\Lambda$).
Two semilattices $S$ and $S'$ in $\Bbb E$ (resp. in $\Lambda$) are said to be
{\it isomorphic} if there exists $\vvp\in GL(\Bbb E)$ 
(resp. $\vvp\in \aut \Lambda$,
the group of automorphisms of $\Lambda$) 
so that $\vvp(S)=S'$,
and denoted $S\cong S'$.
Semilattices $S$ and $S'$ in $\Bbb E$ are said to be
{\it similar} if there exists $\vvp\in GL(\Bbb E)$ 
(resp. $\vvp\in \aut \Lambda$) so that $\vvp(S+\sigma)=S'$
for some $\sigma\in S$,
and we then write $S\sim S'$.
The relations $\cong$ and $\sim$ are equivalence relations.

\example{Example 3.1}
Let $F_{\bs\varepsilon}=\oplus_{\bs\alpha\in\Lambda}\ Ft_{\bs\alpha}$
be an elementary quantum torus. 
We fix
a toral $\langle\bs\sigma_{1},\ldots,\bs\sigma_{n}\rangle$-grading
of $F_{\bs\varepsilon}$.
Let $\tau$ 
be a graded involution of $F_{\bs\varepsilon}$,
and let
$$S(\bs\varepsilon,\tau)
:=\{\bs\alpha\in\Lambda \ |\ \tau(t_{\bs\alpha})=t_{\bs\alpha}\}.$$
Then $S(\bs\varepsilon,\tau)$ satisfies (S1) and (S2),
and so $S(\bs\varepsilon,\tau)$ is a semilattice in some $\Bbb E$.
In [1], p.83, there is a description
of $S(\bs\varepsilon,\tau)$ 
in terms of the coordinates
of $\Lambda$ relative to
the basis
$\langle\bs\sigma_{1},\ldots,\bs\sigma_{n}\rangle$,
namely,
for 
$\bs\alpha
=\alpha_1\bs\sigma_1+\cdots +\alpha_n\bs\sigma_n\in\Lambda$,
$\bs\varepsilon=(\varepsilon_{ij})$ and
$\tau=(a_1,\ldots,a_n)$,
$$S(\bs\varepsilon,\tau)
=\{\bs\alpha\in\Lambda \ |\ 
\sum_{i\in I_{\tau}}\ \alpha_i +\sum_{(i,j)\in J_{\bs\varepsilon}}\ \alpha_i\alpha_j\equiv 0
\mod 2\}$$
where $I_{\tau}=\{ i\ |\ a_i=-1\}$ 
and
$J_{\bs\varepsilon}=\{ (i,j)\ |\ \varepsilon_{ij}=-1\}$.

Now, if $S(\bs\varepsilon,\tau)$ satisfies (S3)$'$,
it is a semilattice in $\Lambda$.
For example,
$S(\bs\varepsilon,*)$ is a semilattice in $\Lambda$
since $\bs\sigma_1,\ldots,\bs\sigma_n\in S(\bs\varepsilon,*)$.
Let
$$\Lambda^{(t)}
=2\Bbb Z\bs\sigma_1+ \cdots +2\Bbb Z\bs\sigma_t+ \Bbb Z\bs\sigma_{t+1}+ \cdots +
\Bbb Z\bs\sigma_n.$$
Then one can see that
$$
S(\bk 1,\tau_1)=\Lambda^{(1)}
\quad\text{and}\quad
S(\bs h_{1,n},\tau_2)=\Lambda^{(2)},
$$
which are lattices, and so semilattices in some Euclidean space 
but not semilattices in $\Lambda$.

If $(F_{\bs\varepsilon},\tau)\cong (F_{\bs\varepsilon'},\tau')$,
then by Lemma 1.2,
there exists the induced automorphism $p$ of $\Lambda$,
and clearly we have
$p\big(S(\bs\varepsilon,\tau)\big)= S(\bs\varepsilon',\tau')$.
Therefore, by Theorem 2.7:
\proclaim{Corollary 3.2}
$$
S(\bs\varepsilon,\tau)\cong
\Lambda^{(1)},\
\Lambda^{(2)}
\quad\text{or}\quad
S(\bs\eta,*)
\quad\text{as semilattices in $\Lambda$}
$$
for some elementary quantum matrix $\bs\eta$.
\endproclaim
\endexample

\enskip

We will need the following fundamental property of
semilattices,
which is shown in [1] II.1.4.

\proclaim{Lemma 3.3}
Suppose that $S$ is a semilattice in a lattice $\Lambda$.
Then 
$$2\Lambda\subset S\subset\Lambda
\quad\text{and}\quad
2\Lambda+S\subset S. \tag3.4$$
Conversely, any generating subset $S$ 
of $\Lambda$ satisfying \rom{(3.4)} is a semilattice in $\Lambda$.
\qed
\endproclaim
Suppose that $S$ is a semilattice in a lattice $\Lambda$.
Then, by (3.4) above, one can write
$$S=\bigsqcup_{i=0}^m\ (\sigma_i+2\Lambda)
\quad \text{(disjoint union)}
\quad
\text{for some $\sigma_i\in S$.}$$
We call the integer $m+1$ the {\it index of $S$} and write it as $I(S)$,
though Azam first defined the index as $m$ (see [4], p.3 Definition 1.5).
We have found our definition more convenient.
Let $n:=\rank \Lambda$.
Then one can check that $n+1\leq I(S)\leq 2^n$.
Azam showed that
the index is a similarity invariant
(see [4] Lemma 1.7, p.3).

\head
\S 4  Classification of
$S(\bs\varepsilon,*)$
\endhead

Recall the notation
$S(\bs\varepsilon,\tau)
=\{\bs\alpha\in\Lambda \ |\ \tau(t_{\bs\alpha})=t_{\bs\alpha}\}$
for a quantum torus
$(F_{\bs\varepsilon},\tau)$
with graded involution,
where $\bs\varepsilon$ is any elementary quantum matrix
and $\tau$ is any graded involution (Example 3.1).
Also, we defined the main involution $*$
of $F_{\bs\varepsilon}$ for any elementary quantum matrix $\bs\varepsilon$,
and two special graded involutions
$\tau_1$ and $\tau_2$ of $F_{\bs h_{l,n}}$ 
for the special elementary quantum matrix $\bs h_{l,n}$
in \S 2.
Note that $n\geq 2l$ and $l\geq 0$. Also, $\tau_1$ is defined when $n>2l$
and $\tau_2$ is defined when $l\geq 1$.

We will classify $S(\bs\varepsilon,*)$ in this section.
By Theorem 2.7, we already know that 
$$S(\bs\varepsilon,*)\cong S(\bs h_{l,n},*),\quad 
S(\bs h_{l,n},\tau_1)\ (l\geq 1) \quad
\text{or}\quad
S(\bs h_{l,n},\tau_2)\ (l\geq 2).$$

For simplicity, we put
$$S(n,l,\tau):=S(\bs h_{l,n},\tau).$$
Let
$F_{\bs h_{l,n}}=F_{\bs h_{l,n}}[t^{\pm 1}_{1},\ldots,t^{\pm 1}_{n}]$
with a $\langle\bs\sigma_{1},\ldots,\bs\sigma_{n}\rangle$-grading.
Let 
$$
\align
I\big(S(n,l,\tau)\big):
&=\{(\kappa_{1},\ldots,\kappa_{n})
\in \{0,1\}^{n}\ |\ 
\kappa_{1}\bs\sigma_{1}+\cdots +\kappa_{n}\bs\sigma_{n}
\in S(n,l,\tau)\} \\
&=\{(\kappa_{1},\ldots,\kappa_{n})
\in \{0,1\}^{n}\ |\ 
\tau(t_{1}^{\kappa_{1}}\cdots t_{n}^{\kappa_{n}})
=t_{1}^{\kappa_{1}}\cdots t_{n}^{\kappa_{n}}\} 
\quad\text{and} \\
I\big(S(n,l,\tau)\big)^{-}:
&=\{0,1\}^{n}\setminus I\big(S(n,l,\tau)\big) \\
&=\{(\kappa_{1},\ldots,\kappa_{n})
\in \{0,1\}^{n}\ |\ 
\tau(t_{1}^{\kappa_{1}}\cdots t_{n}^{\kappa_{n}})
=-t_{1}^{\kappa_{1}}\cdots t_{n}^{\kappa_{n}}\}.
\endalign
$$
So
$$
2^n=| \{0,1\}^{n}|
=| I\big(S(n,l,\tau)\big)|
+| I\big(S(n,l,\tau)\big)^{-}|.
\tag0
$$
We note that
$| I\big(S(n,l,\tau)\big)|$
is the index of
the semilattice $S(n,l,\tau)$ in $\Lambda$ 
if $S(n,l,\tau)=S(n,l,*)$,
$S(n,l,\tau_1)$ for $l\geq 1$
or $S(n,l,\tau_2)$ for $l\geq 2$.
Thus, 
if $| I\big(S(n,l_0,*)\big)|$,
$| I\big(S(n,l_1,\tau_{1})\big)|$
and
$| I\big(S(n,l_2,\tau_{2})\big)|$
are all distinct for any $l_0,l_1,l_2$,
then
the $S(n,l,*)$,
$S(n,l,\tau_1)$
and $S(n,l,\tau_2)$
are pairwise non-similar.
In fact, we can prove the following:
\proclaim{Lemma 4.1}
In the notation above, we have the index formulas
$$
\align
| I\big(S(n,l,*)\big)|
&=2^{n-1}+2^{n-l-1}\quad (l\geq 0),
\\
| I\big(S(n,l,\tau_{1})\big)|
&=2^{n-1} \quad (l\geq 0\ \text{and}\ n>2l)
\\
| I\big(S(n,l,\tau_{2})\big)|
&=2^{n-1}-2^{n-l-1}\quad (l\geq 1).
\endalign
$$
In particular,
for arbitrary $l_0,l_1\geq 0$
and $l_2\geq 1$ such that $n\geq 2l_0,2l_2$
and 
$n>2l_1$,
$$
| I\big(S(n,l_0,*)\big)|
>
| I\big(S(n,l_1,\tau_{1})\big)|
>
| I\big(S(n,l_2,\tau_{2})\big)|.
$$
\endproclaim
\demo{Proof}
For $\bs\kappa=(\kappa_{1},\ldots,\kappa_{n})\in \{0,1\}^{n}$
and
$t^{\bs\kappa}:=t_1^{\kappa_{1}}\cdots t_{2l}^{\kappa_{2l}}
t_{2l+1}^{\kappa_{2l+1}}\cdots t_n^{\kappa_{n}}$,
we have
$$
(t^{\bs\kappa})^*
=(t_2^{\kappa_{2}}t_1^{\kappa_{1}})
(t_4^{\kappa_{4}}t_3^{\kappa_{3}})\cdots 
(t_{2l}^{\kappa_{2l}}t_{2l-1}^{\kappa_{2l-1}})
t_{2l+1}^{\kappa_{2l+1}}\cdots t_n^{\kappa_{n}}
=
(-1)^{\sum_{i=1}^l\ \kappa_{2i-1}\kappa_{2i}}
t^{\bs\kappa}.
$$
Note that
$$
t_{2i}^{\kappa_{2i}}t_{2i-1}^{\kappa_{2i-1}}=
\cases
t_{2i-1}^{\kappa_{2i-1}}t_{2i}^{\kappa_{2i}}&
\text{if $(\kappa_{2i-1},\kappa_{2i})=(0,0)$, $(0,1)$ or $(1,0)$} \\
-t_{2i-1}^{\kappa_{2i-1}}t_{2i}^{\kappa_{2i}}&
\text{if $(\kappa_{2i-1},\kappa_{2i})= (1,1)$}.
\endcases
$$
Hence,
for
$$
\bar l=
\cases
l-1&
\text{if $l$ is even} \\
l&
\text{if $l$ is odd},
\endcases
$$
we obtain,
by counting the pairs
$(\kappa_{2i-1},\kappa_{2i})= (1,1)$,
$$
\align
| I\big(S(n,l,*)\big)|
&=2^n-2^{n-2l}\bigg(
\binom l1 3^{l-1}+\binom l3 3^{l-3}+
\cdots +\binom l{\bar l} 3^{l-\bar l}\bigg)\\
&=2^n-2^{n-2l}(2^{2l-1}-2^{l-1}) \\
&=2^{n-1}+2^{n-l-1},
\tag1
\endalign
$$
by comparing the binomial expansions of $(3+1)^l$ and $(3-1)^l$.

Next we show 
$| I\big(S(n,l,\tau_{1})\big)|=2^{n-1}$
for any $l\geq 0$.
Let
$A_0
:=\{\bs\kappa\in \{0,1\}^{n}\ |\ \kappa_{2l+1}=0\}$ and
$A_1:=\{\bs\kappa\in \{0,1\}^{n}\ |\ \kappa_{2l+1}=1\}$
so that
$$
I\big(S(n,l,\tau_{1})\big)
=\big(I\big(S(n,l,\tau_{1})\big)\cap A_0\big)
\sqcup
\big(I\big(S(n,l,\tau_{1})\big)\cap A_1\big).
$$
Since $\tau_{1}(t_{2l+1})=-t_{2l+1}$
and $t_{2l+1}$
commutes with
all $t_{i}$,
we have
$| \big(I\big(S(n,l,\tau_{1})\big)\cap A_1\big)|
=| I\big(S(n-1,l,*)\big)|$
and
$| \big(I\big(S(n,l,\tau_{1})\big)\cap A_1\big)|
=| I\big(S(n-1,l,*)\big)^{-}|$.
Thus, by (0), we get
$$
| I\big(S(n,l,\tau_{1})\big)|
=| I\big(S(n-1,l,*)\big)|
+| I\big(S(n-1,l,*)\big)^{-}|
=2^{n-1}.
$$

Recall that $\tau_2$ is defined only for $l\geq 1$, and so
 we can consider a partition of $\{0,1\}^{n}$ by
the following four subsets 
$B_k$, $k=1,2,3,4$,
namely,
$$
\align
B_1:
&=\{\bs\kappa\in \{0,1\}^{n}\ |\ \kappa_{2l-1}=\kappa_{2l}=0\}, 
\quad
B_2:=\{\bs\kappa\in \{0,1\}^{n}\ |\ \kappa_{2l-1}=1,\kappa_{2l}=0\}, \\
B_3:
&=
\{\bs\kappa\in \{0,1\}^{n}\ |\ \kappa_{2l-1}=0,\kappa_{2l}=1\}, 
\quad
B_4:=\{\bs\kappa\in \{0,1\}^{n}\ |\ \kappa_{2l-1}=\kappa_{2l}=1\},
\endalign
$$
so that
$$
I\big(S(n,l,\tau_2)\big)
=\bigsqcup_{k=1}^{4}\
\big(I\big(S(n,l,\tau_2)\big)\cap B_k\big).
$$
Since $\tau_{2}(t_{2l-1})=-t_{2l-1}$,
$\tau_{2}(t_{2l})=-t_{2l}$
and $\tau_{2}(t_{2l-1}t_{2l})=-t_{2l-1}t_{2l}$,
and since
$t_{2l-1}$, $t_{2l}$
and $t_{2l-1}t_{2l}$
commute with
all $t_{i}$
for $i\neq 2l-1,2l$,
we have
$| \big(I\big(S(n,l,\tau_2)\big)\cap B_1\big)|
=| I\big(S(n-2,l-1,*)\big)|$
and
$| \big(I\big(S(n,l,\tau_2)\big)\cap B_k\big)|
=| I\big(S(n-2,l-1,*)\big)^{-}|$
for $k=2,3,4$.
Thus we get 
$$
\align
| I\big(S(n,l,\tau_{2})\big)|
&=| I\big(S(n-2,l-1,*)\big)|+3| I\big(S(n-2,l-1,*)\big)^{-}| \\
&=2^{n-2}+2| I\big(S(n-2,l-1,*)\big)^{-}|
\quad\text{by (0)} \\
&=2^{n-2}+2\big(2^{n-2}-(2^{(n-2)-1}+2^{(n-2)-(l-1)-1})\big)
\quad\text{by (0) and (1)} \\
&=2^{n-1}-2^{n-l-1}.
\qed
\endalign
$$
\enddemo
Thus, by the inequalities in Lemma 4.1,
the three semilattices
$$\text{$S(n,l,*)$,
$S(n,l,\tau_1)$ ($l\geq 1$)
and $S(n,l,\tau_2)$ ($l\geq 2$)
are pairwise non-similar in $\Lambda$.}
$$
{\bf End of proof of Theorem 2.7}:
If $(F_{\bs h_{l,n}},*)\cong (F_{\bs h_{l,n}},\tau_{2})$,
then $S(n,l,*)\cong S(n,l,\tau_2)$
as semilattices in $\Lambda$.
Hence as a corollary of Lemma 4.1,
we get
$(F_{\bs h_{l,n}},*)\not\cong (F_{\bs h_{l,n}},\tau_{2})$ 
for $l\geq 2$.
That is, we get
one of the assertions in Theorem 2.7
whose proof was postponed there. \qed

\enskip

Moreover, by the index formulas in Lemma 4.1,
$$\text{$l$ is a similarity invariant for the semilattices $S(n,l,*)$
and $S(n,l,\tau_2)$ ($l\geq 2$) in $\Lambda$.}
$$
To show that $l$ is a similarity invariant
for $S(n,l,\tau_1)$,
we would like to have a new similarity invariant
since the index of $S(n,l,\tau_1)$ is constant for $l\geq 1$.
Thus we define the following:
\definition{Definition 4.2}
Let $S$ be a semilattice in a lattice $\Lambda$.
For $\gamma\in S$, if $\gamma+\sigma\in S$ for all $\sigma\in S$,
then $\gamma$ is called a {\it saturated element of $S$}.
We denote the subset of saturated elements of $S$ by
$\Sigma(S)$.
Then $\Sigma(S)$ is a subgroup of $\Lambda$ 
containing $2\Lambda$.
We define the {\it saturation number $\frak s=\frak s(S)$ of $S$} as
$$|\Lambda/\Sigma(S)|=2^{\frak s}.$$
\enddefinition
\proclaim{Lemma 4.3}
\rom{(i)}
$\Sigma(S)=\Sigma(S+\sigma)$ for 
any semilattice $S$ in $\Lambda$
and any $\sigma\in S$.

\rom{(ii)}
The saturation number is a similarity invariant.
\endproclaim
\demo{Proof}
(i):
Let $\gamma\in \Sigma(S)$.
Then $\gamma-\sigma\in S$ for any $\sigma\in S$, and so $\Sigma(S)\subset S+\sigma$.
Moreover, for the semilattice $S+\sigma$
and any $\rho+\sigma\in S+\sigma$,
we have $\gamma+\rho+\sigma\in S+\sigma$ since $\gamma+\rho\in S$.
Hence $\Sigma(S)\subset \Sigma(S+\sigma)$ for any $\sigma\in S$.
Since $-2\sigma\in S$, we have $-\sigma\in S+\sigma$.
Hence $\Sigma(S+\sigma)\subset \Sigma(S)$,
which shows (i).

(ii): By (i),
we have $\frak s(S)=\frak s(S+\sigma)$ for any $\sigma\in S$.
Hence we only need to show that
the saturation number is an isomorphism invariant.
Suppose $p(S)=S'$ for some $p\in\aut\Lambda$.
Then one can easily see that $p\big(\Sigma(S)\big)=\Sigma(S')$.
Therefore,
$|\Lambda/\Sigma(S)|=|\Lambda/p\big(\Sigma(S)\big)|=|\Lambda/\Sigma(S')|$,
i.e.,
$\frak s$ is an isomorphism invariant.
\qed
\enddemo
\remark{Remark}
One can easily show that
$\Sigma(S)=\bigcap_{\sigma\in S}\ (S+\sigma)$.
\endremark
\proclaim{Corollary 4.4}
Let $l\geq 1$.  Then 
$\Sigma\big(S(n,l,\tau_1)\big)=\Lambda^{(2l+1)}$, and hence
$l$ is a similarity invariant
for the semilattices $S(n,l,\tau_1)$ in $\Lambda$.
\endproclaim
\demo{Proof}
Recall our notation
$S(n,l,\tau_1)=S(\bs h_{l,n},\tau_1)
=\{\bs\alpha\in\Lambda \ |\ \tau_1(t_{\bs\alpha})=t_{\bs\alpha}\}$
for the quantum torus
$F_{\bs h_{l,n}}=F_{\bs h_{l,n}}\lpt$
with
a $\langle\bs\sigma_{1},\ldots,\bs\sigma_{n}\rangle$-grading.
By Lemma 2.5, 
the grading group of the centre $Z(F_{\bs h_{l,n}},\tau_1)$
is equal to $\Lambda^{(2l+1)}$.
Thus it is clear from this that
$\Sigma\big(S(n,l,\tau_1)\big)\supset\Lambda^{(2l+1)}$.
For the other inclusion,
suppose $\Sigma\big(S(n,l,\tau_1)\big)\setminus \Lambda^{(2l+1)}\neq\emptyset$.
Then there exists 
$\bs\kappa:=
\kappa_1\bs\sigma_1+\cdots +\kappa_{2l+1}\bs\sigma_{2l+1}\in \Sigma\big(S(n,l,\tau_1)\big)$,
where $\kappa_i=0$ or $1$
but not all $\kappa_1,\ldots,\kappa_{2l}$ are $0$.
Then 
for $\kappa_j\neq 0$ with $j\leq 2l$,
we have $\bs\sigma_{k}\in S(n,l,\tau_1)$
where
$$k=
\cases
j+1&\text{if $j$ is odd} \\
j-1&\text{if $j$ is even},
\endcases
$$
and $\bs\kappa+\bs\sigma_{k}\notin S(n,l,\tau_1)$
since $\tau_1(t_{1}^{\kappa_1}\cdots t_{2l+1}^{\kappa_{2l+1}}t_k)
=t_kt_{1}^{\kappa_1}\cdots t_{2l+1}^{\kappa_{2l+1}}
=-t_{1}^{\kappa_1}\cdots t_{2l+1}^{\kappa_{2l+1}}t_k$.
This is a contradiction.
Hence $\Sigma\big(S(n,l,\tau_1)\big)=\Lambda^{(2l+1)}$.
Thus $\frak s\big(S(n,l,\tau_1)\big)=2l+1$, and hence
$l$ is a similarity invariant by Lemma 4.3.
\qed
\enddemo
\remark{Remark 4.5}
(i)
One can also check that
$\Sigma\big(S(n,l,*)\big)=\Sigma\big(S(n,l,\tau_2)\big)=\Lambda^{(2l)}$.
So this is another reason why $l$ is a similarity invariant for
$S(n,l,*)$ or $S(n,l,\tau_2)$.

(ii) $S(n,l,\tau_1)$ for $l\geq 1$ give us $[\frac{n}{2}]$ semilattices 
in $\Lambda$ which have the same
index but are not similar,
where $[\frac{n}{2}]$ is the greatest integer less than or equal to $\frac{n}{2}$.
\endremark

\enskip

We summarize the results about the semilattices above as a theorem.
\proclaim{Theorem 4.6}
Let $S(\bs\varepsilon,*)$ be the semilattice in $\Lambda$ defined in Example 3.1.
Then $S(\bs\varepsilon,*)$ is isomorphic to
$$S(\bs h_{l,n},*)\ (l\geq 0),\quad 
S(\bs h_{l,n},\tau_1)\ (l\geq 1) \quad
\text{or}\quad
S(\bs h_{l,n},\tau_2)\ (l\geq 2),$$
and any two of these three semilattices are not similar.
Moreover, 
for each of these
$l$ is a similarity invariant.

In particular, the number of similarity classes of $S(\bs\varepsilon,*)$
is 
$$\cases
3\big[\frac{n}{2}\big]&\text{if $n\geq 4$} \\
2&\text{if $n=2,3$} \\
1&\text{if $n=1$}.
\endcases
$$
\endproclaim
\demo{Proof}
We only need to show the last statement.
Since $l\leq \big[\frac{n}{2}\big]$,
there are $\big[\frac{n}{2}\big]+1$ similarity classes from $S(\bs h_{l,n},*)$
for $n\geq 1$,
$\big[\frac{n}{2}\big]$ classes from $S(\bs h_{l,n},\tau_1)$
for $n\geq 2$
and
$\big[\frac{n}{2}\big]-1$ classes from $S(\bs h_{l,n},\tau_2)$
for $n\geq 4$.
Summing them up, we get the results.
\qed
\enddemo
\remark{Remark 4.7}
The number of similarity classes of semilattices in $\Lambda$
is at least $2^n-n$,
which is bigger than the number above if $n\geq 3$.
Thus if $n$ is not too small, one can say that
the semilattices
$S(\bs\varepsilon,*)$ are far from exhausting all semilattices in $\Lambda$.
\endremark

\head
\S 5 Extended affine root systems of type $\text C$
\endhead

We review the description of extended affine root systems of
type $\text C_r$ for $r\geq 3$ following [1], p.34. 
Let $\Lambda$ be a lattice and $S$ be a semilattice in a Euclidean space $\Bbb E$
so that
$$S+2\Lambda\subset S\quad\text{and}\quad \Lambda+S\subset \Lambda.
\tag5.1$$
Then
an extended affine root system $R$ of type 
$\text C_r$ $(r\geq 3)$
contains an irreducible root system 
$\Delta=\Delta_{sh}\sqcup\Delta_{lg}$ of type $\text C_r$,
where $\Delta_{sh}$ (resp. $\Delta_{lg}$)
is the set of short (resp. long) roots,
so that
$$R=R(\Lambda,S)=\Lambda\sqcup\big(\bigsqcup_{\mu\in\Delta_{sh}}\ (\mu+\Lambda)\big)
\sqcup\big(\bigsqcup_{\mu\in\Delta_{lg}}\ (\mu+S)\big).
\tag5.2
$$
The rank of the lattice $\Lambda$ is called the {\it nullity} of $R$.

If $(\Lambda,S)$ and $(\Lambda',S')$ are pairs of 
a lattice and a semilattice in $\Bbb E$ satisfying (5.1),
we say that $(\Lambda,S)$ and $(\Lambda',S')$ are {\it isomorphic}, written 
$(\Lambda,S)\cong (\Lambda',S')$, if there exists $\vvp\in GL(\Bbb E)$ such that
$\vvp(\Lambda)=\Lambda'$ and $\vvp(S)=S'$.
Also, we say that $(\Lambda,S)$ and $(\Lambda',S')$ are {\it similar}, written 
$(\Lambda,S)\sim (\Lambda',S')$, if there exists $\lambda\in S$
such that
$(\Lambda,S+\lambda)\cong (\Lambda',S')$.
Note that $(\Lambda,S+\lambda)$ is a pair of a lattice and a semilattice satisfying (5.1)
(see Definition 4.8 in [1], p.45).
The relations $\cong$ and $\sim$ are equivalence relations.
It is shown
in [1] Theorem 3.1, p.39
that
the root systems $R(\Lambda,S)$ and $R(\Lambda',S')$
are isomorphic if and only if
$(\Lambda,S)\sim (\Lambda',S')$.

In general, (5.1) implies that $2\Lambda\subset S\subset \Lambda$, and so
$2 \Lambda\subset \langle S\rangle\subset \Lambda$.
Thus we have 
$$| \Lambda/\langle S\rangle|=2^t,
\quad\text{where $0\leq t\leq n$}.$$
The integer $t=t(\Lambda,S)$ is called the {\it twist number}
of the pair $(\Lambda,S)$.
The twist number is a similarity invariant of the pair
(see Definition 4.11 in [1], p.46),
and 
so the twist number
is an isomorphism invariant of the root system
$R(\Lambda,S)$.

\example{Example 5.3}
Let $\Lambda$ be a lattice with basis
$\{\bs\sigma_1,\ldots,\bs\sigma_n\}$.
Then the pair $(\Lambda,\Lambda^{(t)})$ satisfies (5.1)
with twist number $t$,
where $\Lambda^{(t)}$ is defined in Example 3.1.
Moreover, for any semilattice $S'$ in 
$\Bbb Z\bs\sigma_{t+1}+\cdots+\Bbb Z\bs\sigma_n$,
the pair $(\Lambda,2\Bbb Z\bs\sigma_{1}+\cdots+2\Bbb Z\bs\sigma_t+S')$ satisfies (5.1)
with twist number $t$
(Proposition 4.17, p.47 in [1]).
\endexample
 
\enskip

The root systems of extended affine Lie algebras are
extended affine root systems. However, it was conjectured in [1] that
an extended affine root system is not necessarily
the root system of an extended affine Lie algebra.
Allison and Gao have shown in [2] that the twist numbers
of root systems of extended affine Lie algebras of type $\text C_r$
($r\geq 3$) do not exceed $3$.
Precisely, they showed that
such a root system $R$ is given by
$$R\big(\Lambda,S(\bs\varepsilon,\tau)\big)\
\quad\text{if $r\geq 4$},
$$
where 
$S(\bs\varepsilon,\tau)$
is the semilattice of 
$(F_{\bs\varepsilon},\tau)$
for any elementary quantum matrix
$\bs\varepsilon$ 
and any graded involution $\tau$ defined in Example 3.1
and $\Lambda$ is a toral grading of $F_{\bs\varepsilon}$.
If $r=3$, then
$$R\big(\Lambda,S(\bs\varepsilon,\tau)\big)\
\quad\text{or}\quad
R(\Lambda,\Lambda^{(3)}),
$$
where 
the second one
comes from
the octonion torus with standard involution
(see List 6.1, p.46 
and Proposition 4.25, p.20 in [2]).
Then they calculated the twist number of
$\big(\Lambda,S(\bs\varepsilon,\tau)\big)$,
and showed that such numbers do not exceed $2$
(see Theorem 6.2(b), p.46 in [2]).
This fact also follows from our Corollary 3.2.
Namely,
we have
$$
(\Lambda,S(\bs\varepsilon,\tau))\cong
(\Lambda,\Lambda^{(1)}),\
(\Lambda,\Lambda^{(2)})
\quad\text{or}\quad
(\Lambda,S(\bs\eta,*))
$$
for some elementary quantum matrix $\bs\eta$,
and so
$$
t\big(\Lambda,S(\bs\eta,*)\big)=0, \quad
t(\Lambda,\Lambda^{(1)})=1 \quad\text{and}\quad
t(\Lambda,\Lambda^{(2)})=2. 
$$
Note that in general,
even if $t=t(\Lambda,S)=1$, $2$ or $3$,
there are many non-isomorphic semilattices $S$ 
with the same twist number if $n$ is not too small,
as we suggested in Example 5.3.
In fact, if $n\geq 5$,
then there are at least two non-isomorphic semilattices $S$
(exactly two if $t=3$).
However, in the pairs arising from root systems of extended affine Lie algebras,
there is only one,
up to isomorphism, in each case,
i.e.,
$\Lambda^{(1)}$ for $t=1$, $\Lambda^{(2)}$ 
for $t=2$ and $\Lambda^{(3)}$ for $t=3$.

As a corollary of Theorem 4.6, we get:

\proclaim{Corollary 5.4}
Let $R=R(\Lambda,S)$
be the root system of an extended affine Lie algebra of type $C_r$
($r\geq 3$). 
Then if $r\geq 4$, $R$ is isomorphic to
$$
R(\Lambda,S(\bs h_{l,n},*))\ (l\geq 0),\quad
R(\Lambda,S(\bs h_{l,n},\tau_1))\ (l\geq 0) \quad\text{or}\quad
R(\Lambda,S(\bs h_{l,n},\tau_2))\ (l\geq 1),
$$
and if $r=3$, $R$ is isomorphic to
$$
R(\Lambda,S(\bs h_{l,n},*))\ (l\geq 0),\quad
R(\Lambda,S(\bs h_{l,n},\tau_1))\ (l\geq 0), \quad
R(\Lambda,S(\bs h_{l,n},\tau_2))\ (l\geq 1)
\quad\text{or}\quad R(\Lambda,\Lambda^{(3)})
$$
Any two of these root systems are not isomorphic.
Moreover, 
for each of these
$l$ is an isomorphic invariant.

In particular, the number of isomorphism classes of $R$
for $r\geq 4$ (resp. $r=3$)
is 
$$\cases
3\big[\frac{n}{2}\big]+2\ \
(3\big[\frac{n}{2}\big]+3)&\text{if $n\geq 4$} \\
4\ \ (5)&\text{if $n=3$} \\
4\ \ (4)&\text{if $n=2$} \\
2\ \ (2)&\text{if $n=1$}.
\endcases
$$
\endproclaim

\enskip

\enskip

Finally, by Remark 4.7, we have:

\proclaim{Corollary 5.5}
Let $r\geq 3$.
Let
$\Cal R_t$ be the set 
of isomorphism classes of root systems of type $\text C_r$
with nullity $n$ and twist number $t$,
and let $\Cal L\Cal R_t$
be the subset of $\Cal R_t$ consisting
of isomorphism classes of the root systems of extended affine Lie algebras of type $\text C_r$
with nullity $n$ and twist number $t$.
Then 
$\Cal L\Cal R_t=\emptyset$ for all $t> 3$.
Moreover, for $t=0$, $1$, $2$ or $3$,
$\Cal L\Cal R_t$ is a proper subset of $\Cal R_t$
if $n\geq 5$.
\endproclaim

\enddocument